\documentstyle{article}

\input amssym.def
\input amssym.tex

\newcommand{\R}{{\Bbb R}}

\newcommand{\ov}{\overline}
\newcommand{\p}{{\partial}}
\newcommand{\im}{{\sqrt{-1}}}
\newcommand{\cp}{{\C P\,\!^2}}
\newcommand{\ocp}{{\ov{\C P}\,\!^2}}
\newcommand{\op}{{\overline{\partial}}}

\newcommand{\Z}{{\Bbb Z}}
\newcommand{\C}{{\Bbb C}}

\newcommand{\al}{{\alpha}}

\newcommand{\be}{{\beta}}

\newcommand{\Om}{{\Omega}}
\newcommand{\om}{{\omega}}

\newcommand{\de}{{\delta}}
\newcommand{\De}{{\Delta}}

\newcommand{\la}{{\lambda}}

\newtheorem{theorem}{Theorem}[section]

\newtheorem{definition}[theorem]{Definition}
\newtheorem{example}[theorem]{Example}
\newtheorem{remark}[theorem]{Remark}

\newtheorem{proposition}[theorem]{Proposition}
\newcommand{\at}{{@}}

\title{K\"{a}hler Geometry of Toric Varieties and Extremal Metrics}
\author{Miguel Abreu\thanks{Supported by NSF grant DMS 9304580 
through the Institute for Advanced
Study. Current address: Departamento de Matem\'atica, Instituto Superior T\'ecnico,
Av. Rovisco Pais, 1096 Lisboa Codex, Portugal. E-mail: mabreu\at math.ist.utl.pt} \\
Institute for Advanced Study \\ \ \\ revised version}

\begin{document}

\maketitle


\begin{abstract}
A (symplectic) toric variety $X$, of real dimension $2n$, is completely
determined by its moment polytope $\De\subset\R^n$. Recently Guillemin
gave an explicit combinatorial way of constructing ``toric'' K\"{a}hler
metrics on $X$, using only data on $\De$. In this paper, differential
geometric properties of these metrics are investigated using
Guillemin's construction. In particular, a nice combinatorial formula
for the scalar curvature $R$ is given, and the Euler-Lagrange condition
for such ``toric'' metrics being extremal (in the sense of Calabi) is
proven to be $R$ being an affine function on $\De\subset\R^n$.
A construction, due to Calabi, of a $1$-parameter family of extremal
K\"{a}hler metrics of non-constant scalar curvature on $\cp\sharp\ocp$
is recast very simply and explicitly using Guillemin's approach.
Finally, we present a curious combinatorial identity for convex
polytopes $\De\subset\R^n$ that follows from the well known relation
between the total integral of the scalar curvature of a K\"{a}hler
metric and the wedge product of the first Chern class of the underlying
complex manifold with a suitable power of the K\"{a}hler class.
\end{abstract}

\section{Introduction}
\label{sec:intro}

Associated to every convex polytope
$\De\subset\R^n$, satisfying some nondegeneracy and integrality conditions,
there is a unique (modulo suitable equivalence relation) closed connected
symplectic manifold $(X_{\De},\om)$ of dimension $2n$ together with a Hamiltonian
action of the $n$-torus $T^n$ such that the image of the corresponding
moment map $\phi$ is precisely $\De$. The triple $(X_{\De},\om,\phi)$ is called a
(symplectic) toric variety. 
It follows from the construction, due to Delzant~\cite{Del}, 
that $X_{\De}$ comes also equipped with an intrinsic $T^n$-invariant
complex structure compatible with the intrinsic symplectic form. In other
words, $X_{\De}$ has an intrinsic ``toric'' K\"{a}hler metric. In~\cite{Gui2}
Guillemin gave an explicit combinatorial formula, in terms of moment data alone,
for the potential of this K\"{a}hler metric as a function on $\De$ (via $\phi$).

It is natural to ask if these canonical toric metrics have any special differential
geometric properties. One possibility would be a relation with extremal K\"{a}hler
metrics, introduced by Calabi in~\cite{Cal} and~\cite{Cal1}. These metrics are
defined as critical points of the square of the $L^2$-norm of the scalar curvature,
considered as a functional on the space of all K\"{a}hler metrics in a fixed
K\"{a}hler class. The extremal Euler-Lagrange equation is equivalent to the gradient 
of the scalar curvature being an holomorphic vector field (see~\cite{Cal}), 
and these metrics are always invariant under a maximal compact subgroup of 
the group of holomorphic transformations of the underlying complex manifold
(see~\cite{Cal1}). This last property is the main reason to look for a relation
between extremal and toric metrics, the later being also invariant under a big
group of holomorphic transformations (the torus $T^n$).

In this paper we look into this question, following Guillemin's approach 
in~\cite{Gui2}. 
It turns out that the scalar curvature $R$ of a toric metric on $X_{\De}$,
as a function on $\De$, has a nice combinatorial formula, while the
extremal condition is $R$ being an affine function on $\De\subset\R^n$.
A simple analysis of some examples leads quickly to the conclusion that,
except for the case of cartesian products of complex projective spaces,
the canonical toric metric does not satisfy the extremal condition.
However Guillemin's construction not only gives an explicit combinatorial
formula for the potential of the canonical toric metric, but can also
be used to explicitly construct other toric K\"{a}hler metrics on
$X_{\De}$ from data on $\De$ alone. One can try in this way to construct
toric metrics satisfying the extremal condition. Although I do not
know of a general procedure to do so, one can analyse by ``hand'' some simple
yet non-trivial examples in real dimension $4$ ($n=2$). In particular a
construction, due to Calabi~\cite{Cal}, of a $1$-parameter family of 
extremal metrics of non-constant scalar curvature on $\cp\sharp\ocp$ 
can be recast very simply and explicitly using Guillemin's construction.

The paper is organized as follows. In \S\ref{sec:toric} we give some
very brief preliminaries on (symplectic) toric varieties. For more details on
this beautiful subject we recommend Guillemin's recent monograph~\cite{Gui1}.
In \S\ref{sec:metric} we review Guillemin's main construction and
result in~\cite{Gui2} (also explained in detail in an Appendix in~\cite{Gui1}).
The combinatorial formula for the
scalar curvature $R$ of a toric metric on $X_{\De}$ and the extremal
Euler-Lagrange condition are derived in \S\ref{sec:formulas}. 
In \S\ref{sec:combformula} we present a 
combinatorial identity that follows from the well known relation between the
total integral of the scalar curvature of a K\"{a}hler metric and the wedge
product of the first Chern class of the underlying complex manifold with a
suitable power of the K\"{a}hler class. Calabi's family of extremal metrics
on $\cp\sharp\ocp$ is discussed in \S\ref{sec:calfam}.

\medskip

\noindent {\bf Acknowledgements:} This work was carried out while I was a
postdoctoral member of the School of Mathematics of the Institute for
Advanced Study. I take this opportunity to thank the staff and faculty of
the Institute for providing a wonderful environment, both for me and my
family.

I also want to thank V.Guillemin for the interest and encouragement given
to this work, and A.Cannas da Silva and S.Simanca for pointing out corrections
to an earlier version of this paper.

\section{Preliminaries on (symplectic) toric varieties}
\label{sec:toric}

Let $(X,\om)$ be a closed connected $2n$-dimensional symplectic manifold,
and let $\tau : T^n\rightarrow {\rm Diff}(X,\om)$ be an effective
Hamiltonian action of the standard $n$-torus. Let $\phi : X\rightarrow\R^n$
be its moment map. The image $\De$ of $\phi$ is a convex polytope, called
the {\em moment polytope}.

Not every convex polytope in $\R^n$ is the moment polytope of some triple
$(X,\om,\phi)$. The following definition characterizes the ones that are.

\begin{definition} \label{def:delzant} 
A convex polytope $\De$ in $\R^n$ is {\em Delzant} if:
\begin{description}
\item[(1)] there are $n$ edges meeting at each vertex $p$;
\item[(2)] the edges meeting at the vertex $p$ are rational, i.e. each edge
is of the form $p + tv_i,\ 0\leq t\leq \infty,\ {\rm where}\ v_i\in\Z^n$;
\item[(3)] the $v_1, \ldots, v_n$ in (2) can be chosen to be a basis of
$\Z^n$.
\end{description}
\end{definition}

In \cite{Del} Delzant associates to every Delzant polytope $\De\subset\R^n$
a closed connected symplectic manifold $X_\De$ of dimension $2n$ together
with a Hamiltonian $T^n$-action such that the image of the corresponding
moment map is precisely $\De$. Moreover, he shows that this is a bijective
correspondence. More precisely, he proves the following:

\begin{theorem} \label{th:delzant}
Let $(X,\om)$ be a compact, connected, $2n$-dimensional
Hamiltonian $T^n$-space, on which the action of $T^n$ is effective with
moment map $\phi : X\rightarrow\R^n$. Then the image $\De$ of $\phi$ is a
Delzant polytope, and $X$ is isomorphic as a Hamiltonian $T^n$-space to
$X_\De$.
\end{theorem}

\section{The intrinsic K\"{a}hler metric on a (symplectic) toric variety}
\label{sec:metric}

It follows from Delzant's construction that $X_\De$ (or, because of 
Theorem~\ref{th:delzant}, any effective hamiltonian $T^n$-space with the same
moment polytope) is equipped with an intrinsic $T^n$-invariant complex
structure compatible with the intrinsic symplectic form. In other words,
$X_\De$ has an intrinsic K\"{a}hler metric. In~\cite{Gui2} Guillemin gives
an explicit formula for this K\"{a}hler metric in terms of moment data alone.
The details are as follows.

Let $(X,\om)$ be an effective Hamiltonian $T^n$-space, with moment map
$\phi : X \rightarrow\R^n$ and moment polytope $\De = \phi(X)\subset\R^n$.
$\De$ can be described by a set of inequalities of the form $\langle y, u_i
\rangle\geq \la_i,\ i=1,\ldots,d$, each $u_i$ being a primitive element of
the lattice $\Z^n\subset\R^n$ and inward-pointing normal to the $i$th 
$(n-1)$-dimensional face of $\De$. Let 
\begin{eqnarray*}
l_i(y) & = & \langle y,u_i\rangle - \la_i, \ i=1,\ldots,d, \\
l_\infty(y) & = & \sum_{i=1}^{d}\langle y, u_i\rangle
\end{eqnarray*}
and $\De^\circ = {\rm interior\ of}\ \De$. Then $y\in\De^\circ$ if and only
if $l_i(y)>0$ for all $i$.

Guillemin's main result in \cite{Gui2} is the following formula for the
restriction of $\om$ to $\phi^{-1}(\De^\circ)$:
\begin{equation} \label{eq:main}
\om = \im\p\op\phi^\ast(\sum_{i=1}^{d} \la_i(\log l_i) + l_\infty )\ .
\end{equation}
The proof of this formula, which will be more relevant to us then the
formula itself, goes as follows. For the intrinsic complex structure,
$\phi^{-1}(\De^\circ)$ is the complex torus $M=\C\,^n/2\pi\im\Z^n$,
and $T^n = \R^n/2\pi\Z^n$ acts on $M$ by the action:
$$
T^n\times M\rightarrow M,\ (t,z)\mapsto z + \im t\ .
$$
If $\om$ is a $T^n$-invariant K\"{a}hler form on $M$, for which the
$T^n$ action is Hamiltonian, then there exists a function $F=F(x),\ x={\rm Re}
\ z$, such that $\om = 2\im\p\op F$ and the moment map $\phi : M
\rightarrow\R^n$ is given by $\phi(z)=\p F/\p x$.

It follows that $\om$ can be written in the form
\begin{equation} \label{eq:kahlerform}
\frac{\im}{2}\sum_{j,k=1}^{n}\frac{\p^2 F}{\p x_j\p x_k} dz_j\wedge d
\overline{z}_k
\end{equation}
and the restriction to $\R^n = ({\rm Re}\ \C\,^n)$ of the K\"{a}hler metric
is the Riemannian metric
$$
\sum_{j,k=1}^{n}\frac{\p^2 F}{\p x_j\p x_k} dx_j dx_k\ .
$$
Guillemin shows that, under the Legendre transform given by the moment map
\begin{equation} \label{eq:legendrey}
y=\frac{\p F}{\p x} = \phi\ ,
\end{equation}
this is the pull-back of the metric
$$
\sum_{j,k=1}^{n}\frac{\p^2 G}{\p y_j\p y_k} dy_j dy_k
$$
on $\De^\circ$, where
\begin{equation} \label{eq:defg}
G = \frac{1}{2}\sum_{k=1}^{d} l_k(y)\log l_k(y)\ .
\end{equation}
Moreover, the inverse of the Legendre transform (\ref{eq:legendrey}) is
\begin{equation} \label{eq:legendrex}
x_i = \frac{\p G}{\p y_i} + a_i,\ i=1,\ldots,n,
\end{equation}
the $a_i$ being constants. This means that, up to a linear factor in $y$, $G$
is the Legendre function dual to $F$, i.e.
\begin{equation} \label{eq:relxy}
F(x)=\sum_{i=1}^{n}y_i\frac{\p G}{\p y_i} - G(y)
\end{equation}
evaluated at $y=\p F/\p x$.

One can now easily deduce (\ref{eq:main}) from (\ref{eq:kahlerform}),
(\ref{eq:defg}) and (\ref{eq:relxy}).

Note also that it follows from (\ref{eq:legendrey}) and (\ref{eq:legendrex})
that the matrix
\begin{equation} \label{eq:hessy}
(G_{ij}) = (\frac{\p^2 G}{\p y_i\p y_j})
\end{equation}
is the inverse of the matrix
\begin{equation} \label{eq:hessx}
(F_{ij}) = (\frac{\p^2 F}{\p x_i\p x_j})\ ,
\end{equation}
at $y=\p F/\p x$.

This construction of Guillemin can also be used to obtain any other toric
K\"{a}hler metric on $X_{\De}$ from data on $\De$ alone. Indeed, one starts
with any $G':\De^{\circ}\rightarrow\R^n$ of the form
\begin{equation} \label{eq:newdefg}
G' = \frac{1}{2}(\sum_{k=1}^{d} l_k(y)\log l_k(y)) + f(y)\ ,
\end{equation}
where $f$ is a smooth function on some open subset of $\R^n$ containing $\De$,
such that the Hessian of $G'$ is positive definite on $\De^{\circ}$. Using the
above Legendre transformation process in reverse, one gets a function
$F'=F'(x)$ on $M=\C^n / 2\pi\sqrt{-1}\Z^n$ that is the K\"{a}hler potential
for a new K\"{a}hler metric $\om'=2\sqrt{-1}\p\op F'$. Because the behaviour
of $\om'$ at infinity is the same as that of $\om$ (we only changed $G$ by a
nonsingular function $f$ on $\De$), this new K\"{a}hler metric compactifies
in the same way to give a metric on $X_{\De}$. We will see in \S\ref{sec:calfam}
an explicit example of this feature of Guillemin's construction.
\begin{remark}
{\rm In the previous description we are looking at $X_{\De}$ as a fixed complex
manifold $(X_{\De},J)$ with an holomorphic $T^n$-action, having two different
invariant K\"{a}hler forms, the canonical $\om$ and the new $\om'$, for which
the action is Hamiltonian. Since, by construction, the moment polytopes for
$(X_{\De},\om)$ and $(X_{\De},\om')$ are the same $\De\subset\R^n$, we know
from Delzant's Theorem~\ref{th:delzant} that there is an equivariant
symplectomorphism $\psi:(X_{\De},\om)\rightarrow (X_{\De},\om')$.
Using $\psi$ we can view $X_{\De}$ as a fixed symplectic manifold
$(X_{\De},\om)$ with a Hamiltonian $T^n$-action, having two different
invariant compatible complex structures, the canonical $J$ and a new
$J'=\psi^{\ast}(J)$.

This two points of view are completely equivalent, the later being perhaps
more natural in terms of symplectic geometry, while the former is standard
in K\"{a}hler geometry and more convenient for computational purposes. }
\end{remark}

\section{Scalar curvature and the extremal condition}
\label{sec:formulas}

Denote by $(G^{ij})$ and $(F^{ij})$ the inverses of the matrices $(G_{ij})$
and $(F_{ij})$ defined by (\ref{eq:hessy}) and (\ref{eq:hessx}),
with $G$ given by the right-hand side of (\ref{eq:newdefg}).
Recall that $G_{ij} = F^{ij}$ at $y=\p F/\p x$. Denote by $\det G$ and
$\det F$ the determinants of the matrices $(G_{ij})$ and $(F_{ij})$.
Throughout this section repeated indices will be summed from $1$ to $n$.

The scalar curvature $R$ of the K\"{a}hler metric (\ref{eq:kahlerform})
is given by
\begin{equation} \label{eq:scalarx}
R = -\frac{1}{2} F^{ij}\frac{\p^2\log(\det F)}{\p x_i\p x_j}
\end{equation}
(this follows from the standard formulas for the scalar curvature of a
K\"{a}hler metric, see for example~\cite{Cal}, using the fact that in our
case the K\"{a}hler potential $F$ only depends on the real part of the
complex coordinate $z$). To express $R$ as a function on the image of the
moment map $\phi = \p F/\p x$, i.e. as a function of $y$, we first note
that under the change of coordinates $y=\p F/\p x$ we have 
(using~(\ref{eq:legendrex}))
\begin{equation} \label{eq:partialy}
F^{ij}\frac{\p}{\p x_j} = \frac{\p^2 G}{\p y_i\p y_j}\frac{\p}{\p x_j}
= \frac{\p x_j}{\p y_i}\frac{\p}{\p x_j} = \frac{\p}{\p y_i}\ .
\end{equation}
From (\ref{eq:scalarx}) using (\ref{eq:partialy}) we then have that
\begin{eqnarray}
2R & = & - F^{ij}\frac{\p^2\log(\det F)}{\p x_i\p x_j}\  =\  -\frac{\p}{\p y_j}
\frac{\p\log(\det F)}{\p x_j} \nonumber\\
 & = & \frac{\p}{\p y_j}\frac{\p\log(\det G)}{\p x_j}\  =\  \frac{\p}{\p y_j}
(\frac{\p\log(\det G)}{\p y_i}\frac{\p y_i}{\p x_j}) \nonumber\\
 & = & \frac{\p}{\p y_j}(\frac{\p\log(\det G)}{\p y_i} F_{ij})\  =\  
\frac{\p}{\p y_j}(G^{ij}\frac{\p\log(\det G)}{\p y_i})\nonumber \\
 & & \ \nonumber \\
\label{eq:scalary1}
\Rightarrow\  2R & = & \frac{\p}{\p y_j}(G^{ij}\frac{\p\log(\det 
G)}{\p y_i}) \ .
\end{eqnarray}

Formula (\ref{eq:scalary1}) can be simplified using the fact that if
$U=(u^{ab})$ and $Y=(y^{ab})$ are symmetric matrices, with $Y$ positive
definite, then
\begin{equation} \label{eq:algebra}
u^{ab}\frac{\p}{\p y^{ab}}(\log(\det Y)) = {\rm trace}(UY^{-1})\ .
\end{equation}
Indeed, from (\ref{eq:scalary1}) using (\ref{eq:algebra}), we then have
\begin{eqnarray*}
2R & = & \frac{\p}{\p y_i}(G^{ij}\frac{\p\log(\det G)}{\p y_j}) \\
 & = & \frac{\p}{\p y_i}(G^{ij}\frac{\p\log(\det G)}{\p G_{ab}}
       \frac{\p G_{ab}}{\p y_j}) \\
 & = & \frac{\p}{\p y_i}(G^{ij}\frac{\p G_{ab}}{\p y_j}G^{ab})\ = \ 
\frac{\p}{\p y_i}(G^{ij}\frac{G_{jb}}{\p y_a}G^{ab}) \\
 & = & - \frac{\p}{\p y_i}(\frac{\p G^{ij}}{\p y_a}G_{jb}G^{ab})\ = \ 
- \frac{\p}{\p y_i}(\frac{\p G^{ij}}{\p y_a}\de_j^a) \\
 & = & - \frac{\p^2 G^{ij}}{\p y_i\p y_j}\ .
\end{eqnarray*}

The Euler-Lagrange equation defining an extremal K\"{a}hler metric is
given by
$$
({\rm grad}\ R)^{(1,0)} = (F^{ij}\frac{\p R}{\p x_j})\frac{\p}{\p z_i}
\equiv {\rm holomorphic\ vector\ field\ .}
$$
Because the coefficients of the $\p /\p z_i$ are real, the only way the
above vector field is holomorphic is if
$$
F^{ij}\frac{\p R}{\p x_j} = \frac{\p R}{\p y_i} \equiv {\rm constant,}
\ i=1,\ldots,n.
$$

We summarize the previous calculations in the following
\begin{theorem} \label{th:formulas}
Let $\De\subset\R^n$ be a Delzant polytope defined by the inequalities
$\langle y,u_i\rangle\geq\la_i,\ i=1,\ldots,d$, each $u_i$ being a primitive
element of the lattice $\Z^n\subset\R^n$ and inward-pointing normal to 
the $i$th $(n-1)$-dimensional face of $\De$. Let $l_i(y)=\langle y,u_i
\rangle - \la_i$, and define $G$  by the right-hand side of (\ref{eq:newdefg})
and $(G_{ij})$ by (\ref{eq:hessy}). Let $F$ be the Legendre function dual to
$G$ and $(G^{ij}) = (G_{ij})^{-1}$. Then
\begin{description}
\item[(i)] the scalar curvature of the K\"{a}hler metric on
$X_\De$ defined by (\ref{eq:kahlerform}) is given by
\begin{equation} \label{eq:scalary}
R = - \frac{1}{2}\sum_{i,j=1}^{n}\frac{\p^2 G^{ij}}{\p y_i\p y_j}\ ;
\end{equation}
\item[(ii)] this K\"{a}hler metric is extremal if and only if
$$
\frac{\p R}{\p y_i} \equiv {\rm constant,}\ i=1,\ldots,n,
$$
i.e., $R$ is an affine function of $y$.
\end{description}
\end{theorem}

\begin{example} \label{ex:plane}{\rm
For $\cp$ with its standard Fubini-Study K\"{a}hler metric, the linear
functions defining the moment polytope (an equilateral right triangle)
are:
$$
l_1(y)=y_1,\ l_2(y)=y_2,\ l_3(y)=1-y_1-y_2\ .
$$
We then have
$$
G = \frac{1}{2}(y_1\log y_1 + y_2\log y_2 + (1-y_1-y_2)\log(1-y_1-y_2))\ ,
$$
$$
(G_{ij}) = \frac{1}{2(1-y_1-y_2)}\left( 
\begin{array}{cc}
\frac{1-y_2}{y_1} & 1 \\
1 & \frac{1-y_1}{y_2}
\end{array}
\right)\ ,
$$
$$
(G^{ij}) = 2\left(
\begin{array}{cc}
y_1(1-y_1) & -y_1 y_2 \\
- y_1 y_2 & y_2(1-y_2)
\end{array}
\right)
$$
and $R = 2 + 1 + 1 + 2 = 6$. }
\end{example}

\section{A combinatorial formula}
\label{sec:combformula}

It is well known that on any K\"{a}hler manifold $X$ the total scalar curvature
is a topological invariant, depending only on the K\"{a}hler class $\Om$ 
represented by the metric $\om$ and the first Chern class $c_1(X)$ of the
underlying complex manifold. More precisely we have that
\begin{equation} \label{eq:totalscalar}
\ov{R} = \int_{X} R_{\om} dv_{\om} = \frac{2\pi}{(n-1)!}\int_{X} c_1(X)
\wedge\om^{n-1} = \frac{2\pi}{(n-1)!} c_1(X)\wedge\Om^{n-1}\ .
\end{equation}

For a toric variety $X_{\De}$, associated to a Delzant polytope $\De\subset\R^n$,
the right-hand side of the above equation can be computed in very simple terms
from combinatorial data on $\De$.  
Recall that $\De$ is defined by the system of inequalities
$\langle y,u_i\rangle\geq\la_i,\ i=1,\ldots,d$, each $u_i$ being a primitive
element of the lattice $\Z^n\subset\R^n$ and inward-pointing normal to 
the $i$th $(n-1)$-dimensional face of $\De$. Let $v_i = -u_i$, the outward-pointing
normal to this face, and let $s_i = -\la_i$. Then $\De$ can also be defined by the
inequalities:
$$
\langle y,v_i\rangle\leq s_i,\ i=1,\ldots,d\ .
$$
Let $v(s_1,\ldots,s_d)$ denote the Euclidean volume of this set.

As before, define $l_i(y)=\langle y,u_i\rangle -\la_i = s_i - \langle y,v_i\rangle,
\ i=1,\ldots,d$, and let $\phi:X_{\De}\rightarrow\De\subset\R^n$ be the moment map.
The subset $X_i\subset X_{\De}$ defined by $l_i\circ\phi = 0$ is a complex
submanifold of real codimension $2$. It is the pre-image in $X$ of the $i$th
$(n-1)$-dimensional face of $\De$. Let $c_i$ be the cohomology class in
$H^2(X_{\De},\Z)$ dual to the homology class $[X_i]$ in $H_{2n-2}(X_{\De},\Z)$.

We then have that (see~\cite{Gui2} and~\cite{Dan})
\begin{equation} \label{eq:kclass}
\frac{\Om}{2\pi} = \frac{[\om]}{2\pi} = \sum_{i=1}^{d} s_i c_i
\end{equation}
and
\begin{equation} \label{eq:chernclass}
c_1(X_{\De}) = \sum_{i=1}^{d} c_i\ .
\end{equation}
Moreover, it follows from (\ref{eq:kahlerform}) and the change of coordinates
(\ref{eq:legendrey}) that
\begin{eqnarray*}
 & & \int_{X_{\De}}\exp(\om) = \int_{X_{\De}}\frac{\om^n}{n!} =
(2\pi)^n {\rm\ Volume }(\De) \\
 & \Rightarrow & v(s_1,\ldots,s_d) = \int_{X_{\De}}\exp(\sum s_i c_i) \\
 & \Rightarrow & \sum_{i=1}^{d}\frac{\p v}{\p s_i} = \int_{X_{\De}}
(\sum c_i)\exp(\sum s_i c_i)\ .
\end{eqnarray*}
Using (\ref{eq:kclass}) and (\ref{eq:chernclass}) we get
$$
\frac{2\pi}{(n-1)!} c_1(X_{\De})\wedge\Om^{n-1} = (2\pi)^n 
\sum_{i=1}^{d}\frac{\p v}{\p s_i}\ ,
$$
the desired combinatorial formula for the right-hand side of
(\ref{eq:totalscalar}).

Regarding its left-hand side, using again (\ref{eq:kahlerform}) and
(\ref{eq:legendrey}), and noting that the scalar curvature $R$, as a
function on $\phi^{-1}(\De^{\circ}) = \C^n / 2\pi\sqrt{-1}\Z^n$, depends
only on $x = {\rm Re\ }z$, we have that
$$
\int_{X_{\De}} R dv = \int_{X_{\De}} R\exp(\om) = (2\pi)^n\int_{\De} R(y) dy\ .
$$

Taking into account (\ref{eq:scalary}), we have proved the following
\begin{proposition} \label{prop:combformula}
Let $\De = \De(s_1,\ldots,s_d)\subset\R^n$ be a Delzant polytope defined by the
inequalities $\langle y,v_i\rangle\leq s_i,\ i=1,\ldots,d$, each $v_i$ being a
primitive element of the lattice $\Z^n\subset\R^n$ and outward-pointing normal
to the $i$th $(n-1)$-dimensional face of $\De$. Let $l_i(y)=s_i - \langle
y,v_i\rangle$, and define $G$  by the right-hand side of (\ref{eq:newdefg})
and $(G_{ij})$ by (\ref{eq:hessy}). Let $(G^{ij}) = (G_{ij})^{-1}$ and
$v(s_1,\ldots,s_d)$ be the
Euclidean volume of $\De(s_1,\ldots,s_d)$. Then
\begin{equation} \label{eq:combformula}
\sum_{i=1}^{d}\frac{\p v}{\p s_i} = - \frac{1}{2}\int_{\De}
(\sum_{i,j}\frac{\p^2 G^{ij}}{\p y_i\p y_j}) dy\ .
\end{equation}
\end{proposition}
It would be interesting to know if (\ref{eq:combformula}) is more than a
curiosity, possibly having some useful application in particular when
$G$ has the canonical form given by (\ref{eq:defg}).

\section{Calabi's family of extremal metrics on ${\bf C}P^2\sharp
\ov{{\bf C}P}^2$}
\label{sec:calfam}

In~\cite{Cal} Calabi constructed families of extremal metrics on
holomorphic $\C P^1$ bundles over $\C P^n$, parametrized by the
cohomology classes of the corresponding K\"{a}hler forms. When these
bundles are holomorphically non-trivial the corresponding group of
complex automorphisms is not reductive, and so, due to a theorem of
Matsushima~\cite{Mat} and Lichn\'erowicz~\cite{Lic}, we know apriori that
these metrics have non-constant scalar curvature. The purpose of
this section is to recast Calabi's simplest example, $\cp\sharp\ocp$,
in the language of the previous sections. For yet another approach to
Calabi's families of extremal metrics see~\cite{Sim}.

The linear functions defining the moment polytope are now:
$$
l_1(y)=y_1,\ l_2(y)=y_2,\ l_3(y)=1-y_1-y_2,\ l_4(y)=y_1+y_2-a,
$$
where $0<a<1$ parametrizes the amount by which we are blowing up the
standard $\cp$. It is easy to check that the metric corresponding, under
Guillemin's construction, to this linear functions is not extremal.
As explained at the end of \S\ref{sec:metric}, in order to get an extremal
metric within the same K\"{a}hler class we can try to change this one
by adding to the corresponding $G$ a smooth function on the
corresponding $\De$ . The fact that the maximal compact subgroup of the group
of complex automorphisms of $\cp\sharp\ocp$ is $U(2)$, and not just $T^2$,
translates into the fact that the smooth function we are looking for will
only depend on the variable $\psi=y_1+y_2$. Hence we will consider $G$ of
the form:
\begin{equation} \label{eq:cpotential}
G=\frac{1}{2}(\sum_{i=1}^{4}l_i(y)\log l_i(y) + f(\psi))\ .
\end{equation}

Calabi's construction is on $\C^2\setminus\{0\}$, say with coordinates
$w=(w_1,w_2)$. His $U(2)$-invariant extremal K\"{a}hler metric is generated
by a K\"{a}hler potential $\Phi(w,\ov{w})=u(t),\ t=\log(|w_1|^2 + |w_2|^2)$,
where $u:\R\rightarrow\R$ is a smooth function satisfying
\begin{equation} \label{eq:nondeg}
u'(t)>0,\ u''(t)>0
\end{equation}
and
\begin{equation} \label{eq:bconditions}
a=\lim_{t\rightarrow -\infty} u'(t),\ 1=\lim_{t\rightarrow +\infty} u'(t)\ .
\end{equation}
The scalar curvature is given by
$$
R=12c_1 u'(t) + 6 c_2
$$
where $c_1$ and $c_2$ are constants given by
$$
c_1 = \frac{2a}{(1-a)(1+4a+a^2)}\ \ {\rm and}\ \ 
c_2 = \frac{1-3a^2}{(1-a)(1+4a+a^2)}\ .
$$
On $\De$ we know from Theorem~\ref{th:formulas} (ii)
that $R$ is an affine function of $\psi=y_1+y_2$, i.e.
$u'(t)=\al\psi+\be$ for some $\al,\be\in\R$. 
Because of the ``boundary conditions'' 
(\ref{eq:bconditions}) we actually have that
\begin{equation} \label{eq:relty}
u'(t)=\psi=y_1+y_2\ .
\end{equation}
Equation (\ref{eq:relty}), together with the non-degeneracy condition
(\ref{eq:nondeg}), gives an implicit formula for $t$ as a function of
$\psi$. The extremal condition in Calabi's construction is then given by
\begin{equation} \label{eq:cextremal}
\frac{dt}{d\psi} = \frac{1-a}{(\psi-a)(1-\psi)} + \frac{2a(1-a)}
{2a\psi^2 + (1+2a-a^2)\psi + 2a^2}\ .
\end{equation}

We are now in position to determine $f(\psi)$. From (\ref{eq:relxy}) and
(\ref{eq:cpotential}) we have that
\begin{eqnarray*}
\frac{1}{2}u(t) & = & \sum_{i=1}^{n}y_i\frac{\p G}{\p y_i} - G(y) \\
 & = & \frac{1}{2}(\psi + \psi f'(\psi) - f(\psi) - \log(1-\psi) +
a\log(\psi-a))\ .
\end{eqnarray*}
Taking the derivative with respect to $\psi$ and using (\ref{eq:relty}) and
(\ref{eq:cextremal}) gives
$$
\psi(\frac{1-a}{(\psi-a)(1-\psi)} + \frac{2a(1-a)}{2a\psi^2 + (1+2a-a^2)\psi 
+ 2a^2}) = 1 + \psi f''(\psi) + \frac{(1-a)\psi}{(1-\psi)(\psi-a)}
$$
which immediately implies
\begin{equation} \label{eq:deff}
f''(\psi) = \frac{2a(1-a)}{2a\psi^2 + (1+2a-a^2)\psi + 2a^2} - \frac{1}{\psi}\ ,
\end{equation}
a formula that can be explicitly integrated twice.
It is now elementary, though tedious, to verify that using
(\ref{eq:cpotential}), (\ref{eq:deff}) and (\ref{eq:scalary}) one indeed gets
$$
R = 12c_1\psi + 6c_2\ .
$$


\newpage

\begin{thebibliography}{99}

\bibitem{Cal}  E.Calabi, Extremal K\"{a}hler metrics, in {\it Seminar on
               Differential Geometry}, ed. S.T.Yau, Annals of Math. Studies 
               {\bf 102}, 159-290, Princeton Univ. Press, 1982.
\bibitem{Cal1} E.Calabi, Extremal K\"{a}hler metrics II, in {\it Differential
               Geometry and Complex Analysis}, eds. I.Chavel and H.M.Farkas,
               Springer-Verlag, 95-114, 1985.
\bibitem{Dan}  V.I.Danilov, The geometry of toric varieties, {\it Russian Math.
               Surveys} {\bf 33} (1978), 97-154.
\bibitem{Del}  T.Delzant, Hamiltoniens p\'{e}riodiques et image convex de
               l'application moment, {\it Bull. Soc. Math. France},
               {\bf 116} (1988), 315-339.
\bibitem{Gui1} V.Guillemin, {\it Moment Maps and Combinatorial Invariants of
               Hamiltonian $T^n$-spaces}, Progress in Math. {\bf 122},
               Birkh\"{a}user, 1994.
\bibitem{Gui2} V.Guillemin, K\"{a}hler structures on toric varieties,
               {\it Journal of Differential Geometry} {\bf 40} (1994), 285-309.
\bibitem{Lic}  A.Lichn\'erowicz, Sur les transformations analytiques des
               vari\'eties k\"{a}hl\'eriennes, {\it C. R. Acad. Sci. Paris}
               {\bf 244} (1957), 3011-3014.
\bibitem{Mat}  Y.Matsushima, Sur la structure du groupe d'hom\'eomorphismes
               analytiques d'une certaine vari\'et\'e k\"{a}hl\'erienne,
               {\it Nagoya Math. Journal} {\bf 11} (1957), 145-150.
\bibitem{Sim}  S.Simanca, A note on extremal metrics of non-constant scalar
               curvature, {\it Israel Journal Math.} {\bf 78} (1992), 85-93.

\end{thebibliography}
\end{document}